\documentclass[a4paper,11pt]{article}
\usepackage[latin1]{inputenc}
\usepackage[T1]{fontenc}
\usepackage{fontenc}
\usepackage{color}
\usepackage{amsmath}
\usepackage{graphicx}
\usepackage{array}
\usepackage{amsfonts}
\usepackage{amssymb}
\usepackage{euscript}
\usepackage{delarray}
\usepackage{latexsym}
\usepackage{enumerate}
\usepackage{float}

\title{Bayesian analysis of hierarchical multi-fidelity codes.}
\author{Loic Le Gratiet  $^\dag$  $^\ddag$ \\ $^\dag$ Universit\'e Paris Diderot 75205 Paris Cedex 13 \\ $^\ddag$ CEA, DAM, DIF, F-91297 Arpajon, France  \\ loic.le-gratiet@cea.fr}

\newcommand{\normal}[2]{\mathcal{N}\left(#1,#2 \right)}

\begin{document}
\maketitle

\section{Abstract}

This paper deals with the Gaussian process based approximation of a code which can be run at different levels of accuracy. This   method, which is a particular case of co-kriging,  allows us to improve a surrogate model of a complex computer code using fast approximations of it. In particular, we focus on the case of a large number of code levels on the one hand and on a Bayesian approach when we have two levels on the other hand.

The main results of this paper are a new approach to estimate the model parameters which provides a closed form expression for an important parameter of the model (the scale factor), a reduction of the numerical complexity  by simplifying the covariance matrix inversion, and a new Bayesian modelling that gives an explicit representation of the joint distribution of the parameters and that is not computationally expensive.

A thermodynamic example is used to illustrate the comparison between 2-level and 3-level co-kriging.
Keywords: surrogate models, co-kriging, multi-fidelity computer experiment, Bayesian analysis.

\section{Introduction}

Large computer codes are widely used in science and engineering to study physical systems since real experiments are often costly and sometimes impossible. Nevertheless, simulations can  sometimes be costly and time-consuming as well. In this case, conception based on an exhaustive exploration of the input space of the code is generally impossible under reasonable time constraints. Therefore, a mathematical approximation of the output of the code - also called surrogate or metamodel - is often built with a few simulations to represent the real system. 

Gaussian Process regression is a particular class of surrogate which makes the assumption that prior beliefs about the code can be modelled by a Gaussian Process. We focus here on this metamodel and on its extension to multiple response models. The reader is refered to \cite{S03} and \cite{R06} for further detail about Gaussian Process models.

Actually, a computer code can often be run at different levels of complexity and a hierarchy of  levels of code can hence be obtained. The aim of this paper is to study the use of several levels of a code to predict the output of a costly computer code.

A first metamodel for multi-level computer codes was built by Kennedy and O'Hagan \cite{KO00} using a spatially stationary correlation structure. This multi-stage model is a particular case of co-kriging which is a well known geostatistical method. Then, Forrester et al. \cite{AF07} went into more detail about the estimation of the model parameters. Furthermore, Forrester \emph{et al.} presented the use of co-kriging for multi-fidelity optimization based on the EGO (Efficient Global Optimization) algorithm created by Jones et al. \cite{DRJ98}. A Bayesian approach was also proposed by Qian and Wu \cite{QW07} which is computationally expensive and does not provide explicit formulas for the joint distribution of the parameters.

This paper presents a new approach to estimate the parameters of the model which is effective in the case of non-spatial stationarity and when many levels of codes are available. In particular, it provides a closed form expression for the estimation of the scale factor which is new and of great practical interest. Furthermore, this approach allows us to consider prior information in the estimation of the parameters. We also address the problem of the inversion of the co-kriging covariance matrix when the number of levels is large. A solution to this problem is provided which shows that the inverse can be easily calculated. Finally, it is known that with a non-Bayesian approach, the variance of the predictive distribution may be underestimated \cite{KO00}. This paper suggests a Bayesian modelling different from the one presented in \cite{QW07}  which provides an explicit representation of the joint distribution for the parameters and avoids prohibitive implementation.

\section{Building a surrogate model based on a hierarchy of $s$ levels of code}\label{building}

Let us assume that we have $s$ levels of code $z_1(x),\dots,z_s(x)$, $x \in \mathbb{R}^d, d > 0$. For all $t=1,\dots,s$ the $t^{th}$ scalar output $z_t(x)$  is modelled by $z_t(x) = Z_t(x,\omega)$ where $Z_t(x,\omega), \omega \in \Omega$ is a realization of the Gaussian process $Z_t(x)$. We will introduce below a consistent set of hypotheses so that the joint process $(Z_t(x))_{x \in \mathbb{R}^d,t=1,\dots,s }$ is Gaussian given a certain set of parameters. Kenney and O'Hagan  \cite{KO00} suggest an autoregressive model to build a metamodel based on a multi-level computer code. Hence, we have a hierarchy of $s$ levels of code - from the less accurate to the most accurate - and for each level, the conditional distribution of the Gaussian process $Z_t(x)$ knowing $Z_1(x),\dots,Z_{t-1}(x)$ is entirely determined by $Z_{t-1}(x)$. Let us introduce here the mathematical formalism that we will use in this paper.

$Q \subset \mathbb{R}^d$ is a  compact subset of $\mathbb{R}^d$ called the input space or the domain of interest.
For $ t = 1,\dots,s$, $D_t =\{x_1^{(t)},\dots,x_{n_t}^{(t)} \}$ is the experimental design set at level $t$ containing $n_t$ points in $Q$. Let $\mathcal{Z}_t = Z_t(D_t) = (Z_t(x_1^{(t)}),\dots,Z_t(x_{n_t}^{(t)}))^T$ be the random Gaussian vector containing the values of $Z_t(x)$ for $x \in D_t$ where  $^T$ stands for the transpose. Let $\mathcal{Z} = (\mathcal{Z}_1^T,\dots,\mathcal{Z}_s^T)^T$ be the Gaussian random vector containing the values of the processes $(Z_t(x))_{t=1,\dots,s}$ at the points of the design sets $(D_t)_{t=1,\dots,s}$. We assume here that the code output is observed without measurement error. The column vector of responses is written $z = (z_1^T,\dots,z_{s}^T)^T$, where $z_t = (z_t(x_1^{(t)}),\dots,z_t(x_{n_t}^{(t)}))^T$ is the output vector for the level $t$.

If we consider $Z_s(x)$, the Gaussian process modelling the most accurate code, we want to determine the predictive distribution of $Z_s(x_0)$, $x_0 \in Q$ given $\mathcal{Z} = z$, \emph{i.e.} the following conditional distribution: $[Z_s(x_0)|\mathcal{Z}= z]$.  

We assume the Markow property introduced in \cite{KO00}:
\begin{equation}
\mathrm{Cov}(Z_t(x),Z_{t-1}(x')|Z_{t-1}(x)) = 0 \qquad \forall x \neq x'
\end{equation}
This means that if $Z_{t-1}(x)$ is known, then nothing more can be learnt about $Z_t(x)$ from any other run of the cheaper code $Z_{t-1}(x')$ for $x' \neq x$.
This assumption leads to the following autoregressive model:
\begin{equation}\label{eq01}
Z_t(x) = \rho_{t-1}(x)Z_{t-1}(x) + \delta_t(x) \quad t = 2, \dots, s
\end{equation}
where $\delta_t(x)$ is a Gaussian process independent of $Z_{t-1}(x),\dots,Z_1(x)$ and $\rho_{t-1}(x)$ represents a scale factor between $Z_t(x)$ and $Z_{t-1}(x)$. We assume that $\rho_{t-1}(x)= f_{\rho_{t-1}}(x)^T\beta_{\rho_{t-1}}$, $ t = 2, \dots, s $,
where $f_{\rho_{t-1}}(x) = (f_{\rho_{t-1}}^1(x),\dots,f_{\rho_{t-1}}^{q_{t-1}}(x))^T$ is a vector of $q_{t-1}$ regression functions - generally including the constant function : $x \in Q \rightarrow 1$ - and $ \beta_{\rho_{t-1}} \in \mathbb{R}^{q_{t-1}}$.\\
Conditioning on parameters $\sigma_t$, $\beta_t$ and $\theta_t$, $\delta_t(x)$ is assumed to be a  Gaussian process with mean $f_t(x)^T \beta_t$, where $f_t(x) $ is a $p_t$-dimensional vector of regression functions, and with a covariance function of the form $c_t(x,x') = \mathrm{cov}(\delta_t(x),\delta_t(x')) = \sigma_t^2r_t(x-x' ; \theta_t)$, where $\sigma_t^2$ is the variance of the Gaussian process and $\theta_t$ are the hyper parameters of the correlation function $r_t$. Moreover, conditioning on parameters $\sigma_1$, $\beta_1$ and $\theta_1$, the simplest code $Z_1(x)$ is modelled as a  Gaussian process with mean $f_1(x)^T \beta_1$ and with covariance function $c_1(x,x') = \sigma_1^2r_1(x-x' ; \theta_1)$. With this consistent set of hypotheses, the joint process $(Z_1(x), \dots, Z_t(x))_{x\in Q, t=1,\dots,s}$ given $\sigma^2 = (\sigma_i^2)_{i=1,\dots,t}$,  $\theta = (\theta_i)_{i=1,\dots,t}$,  $\beta = (\beta_i)_{i=1,\dots,t}$ and $\beta_\rho = (\beta_{\rho_{i-1}})_{i=2,\dots,t}$,  is Gaussian with mean:

\begin{equation}\label{hsx}
\mathbb{E}[Z_t(x)|\sigma^2,\theta,\beta,\beta_\rho] = h_t'(x)^T \beta
\end{equation}

\begin{equation}
h_t'(x)^T = \left( \left(\prod_{i=1}^{t-1} {\rho_{i}(x) }\right)f_1^T(x) , \left(\prod_{i=2}^{t-1} {\rho_{i}(x) }\right)f_2^T(x) , \dots , \rho_{t-1}(x) f_{t-1}^T(x) , f_t^T(x) \right)
\end{equation}

and covariance:

\begin{equation}\label{s2sx}
\mathrm{cov}(Z_t(x),Z_t(x')|\sigma^2,\theta,\beta,\beta_\rho) =
\sum_{j=1}^{t}{\sigma_{j}^2\left( \prod_{i=j}^{t-1} {\rho_{i}^2(x) } \right)r_j(x-x';\theta_j)}
\end{equation}

For each level $t=2,\dots,s$, the experimental design $D_t$ is assumed to be such that $D_t \subseteq D_{t-1}$. Note that this assumption is not necessary but allows us to have closed form expressions for the parameter estimation formulas. Furthermore, we denote by $R_t(D_k,D_l)$ the correlation matrix between points in $D_k$ and $D_l$, $1 \leq k,l \leq s$. $R_t(D_k,D_l)$ is a $(n_k \times n_l)$ matrix with $(i,j)$ entry given by:
\begin{displaymath}
[R_t(D_k,D_l)]_{i,j} = r_t(x_i^{(k)}-x_j^{(l)} ; \theta_t) \qquad 1 \leq i \leq n_k \quad 1 \leq j \leq n_l  
\end{displaymath}
We will use the notation: $R_t(D_k) = R_t(D_k,D_k) $.\\
\cite{KO00}  present the case where $\forall t\in [2,s]$, $\rho_{t-1}(x) = \rho_{t-1}$ is constant. Here, we will consider the general model presented in equations (\ref{eq01}). We will also propose a new approach to estimate the coefficients $(\beta_t,\beta_{\rho_{t-1}})_{t=2,\dots,s}$ based on a Bayesian approach, which allows us to get information about their uncertainties. In the following section, we describe the case of 2 levels of code where the scaling coefficient $\rho$ is constant and then we will extend it for $s$ levels in Section \ref{CaseSLevels}. The general case in which $\rho$ depends on $x$  is addressed in Appendix \ref{appendix1}.

\section{Building a model with 2 levels of code}

Let us assume that we have 2 levels of code $z_2(x)$ and $z_1(x)$. From the previous section we  assume that:
\begin{equation}\label{eq03}
\left\lbrace \begin{array}{c}
Z_2(x)=\rho Z_1(x)+\delta(x), \quad x\in Q \\
(Z_1(x))_{x\in Q}\perp (\delta(x))_{x\in Q} \\
\end{array}\right. 
\end{equation}

The goal of this section is to build a surrogate model for $Z_2(x)$ given the observations $\mathcal{Z}=z$ with an uncertainty quantification. The strategy is the following one. In Subsection \ref{CondDistOut} we describe  the statistical distribution of the output $Z_2(x_0)$ at a new point $x_0$ given the parameters $ (\beta_1,\beta_2,\rho),$ $(\sigma_1^2,\sigma_2^2)$ and $ (\theta_1,\theta_2)$ and the observations $z$. In Subsection  \ref{Bayesian_estimation_of_parameters} we describe the Bayesian estimation of the parameters $ (\beta_1,\beta_2,\rho)$ and $(\sigma_1^2,\sigma_2^2)$ given the observations. As pointed out at the end of Subsection \ref{Bayesian_estimation_of_parameters} the hyper-parameters $(\theta_1,\theta_2)$ are estimated using a concentrated restricted log-likelihood method.

\subsection{Conditional distribution of the output}\label{CondDistOut}

For a point $x_0 \in Q$ we  determine in this subsection the distribution of $[Z_2(x_0)|\mathcal{Z} = z, (\beta_1,\beta_2,\rho),$ $(\sigma_1^2,\sigma_2^2), (\theta_1,\theta_2)]$. Standard results for normal distributions give that:
\begin{equation}\label{eq03bis}
[Z_2(x_0)|\mathcal{Z} = z, (\beta_1,\beta_2,\rho),(\sigma_1^2,\sigma_2^2), (\theta_1,\theta_2)] \sim \mathcal{N}(m_{Z_2}(x_0),s_{Z_2}^2(x_0))
\end{equation}
with mean function:
\begin{equation}
m_{Z_2}(x) = h'(x)^T\beta + t(x)^TV^{-1}(z-H\beta)
\end{equation}
and variance:
\begin{equation}\label{eq37}
s_{Z_2}^2(x)= \rho^2 \sigma_1^2 + \sigma_2^2 - t(x)^TV^{-1}t(x)
\end{equation}
where we have denoted $
\beta = 
\left(
\begin{array}{c}
\beta_1 \\
\beta_2 \\
\end{array}
\right) 
$, $
z = 
\left(
\begin{array}{c}
z_1 \\
z_2 \\
\end{array}
\right) 
$
and where $H$ is defined by:
\begin{displaymath}\label{equationH}
H=\left( \begin{array}{cc}
f_1^T(x_1^{(1)})& 0 \\
\vdots & \vdots\\
f_1^T(x_{n_1}^{(1)}) & 0 \\
\rho f_1^T(x_1^{(2)}) & f_2^T(x_1^{(2)})\\
\vdots & \vdots \\
\rho f_1^T(x_{n_2}^{(2)})& f_2^T(x_{n_2}^{(2)})\\
\end{array}\right) 
=
\left( \begin{array}{c|c}
& \\
F_1(D_1)& 0 \\
& \\
\hline
& \\
\rho F_1(D_2) & F_2(D_2)\\
& \\
\end{array}\right) 
\end{displaymath}
with the notation 
$F_i(D_j) = \left(
\begin{array}{c}
f_i^T(x_{n_1}^{(j)}) \\
\vdots \\
f_i^T(x_{n_j}^{(j)})
\end{array}
\right)$.
 Furthermore, we have $h'(x) = \left(\rho f_1^T(x),f_2^T(x) \right)$  and:
\begin{equation}
\begin{array}{rl}
t(x)^T & = \mathrm{Cov}(Z_2(x),\mathcal{Z}) \\ & = \left( \rho \sigma_1^2 R_1(\{x\},D_1) , \rho^2  \sigma_1^2 R_1(\{x\},D_2)+\sigma_2^2 R_2(\{x\},D_2) \right) \\
\end{array}
\end{equation}
The covariance matrix $V$ of the Gaussian vector $ \mathcal{Z} = \left(\begin{array}{c}\mathcal{Z}_1 \\ \mathcal{Z}_2 \\\end{array}\right)$ can be written :
\begin{equation}
V=\left( 
\begin{array}{cc}
\sigma_1^2 R_1(D_1) & \rho \sigma_1^2  R_1(D_1,D_2) \\
\rho \sigma_1^2  R_1(D_2,D_1) & \rho^2 \sigma_1^2 R_1(D_2)+\sigma_2^2R_2(D_2) \\
\end{array}
\right) 
\end{equation}

\subsection{Bayesian estimation of the parameters with 2 levels of code}\label{Bayesian_estimation_of_parameters}

In this Subsection, we describe the estimation of the parameters $(\beta_1,\beta_2,\rho,\sigma_1^2,$ $\sigma_2^2,\theta_1,\theta_2)$ for the 2-level model given the observations $\mathcal{Z}=z$. Due to the conditional independence  between $Z_1(x)$ and $\delta(x)$, it is possible  to estimate separately the parameters $(\beta_1,\sigma_1^2,\theta_1)$ and $(\beta_2,\rho,\sigma_2^2,\theta_2)$. We first describe the estimations of $(\beta_1,\sigma_1^2)$ given $\theta_1$ and  $(\beta_2,\sigma_2^2, \rho)$ given $\theta_2$, which can be obtained in closed forms. We then describe how to estimate $\theta_1$ and $\theta_2$.

Firstly, we consider the parameters $(\beta_1,\sigma_1^2,\theta_1)$. We  choose the following non-informative prior distributions corresponding to the ``Jeffreys priors" \cite{J61}:
\begin{equation}\label{JeffPrior01}
p(\beta_1|\sigma_1^2,\theta_1) \propto 1 \qquad p(\sigma_1^2) \propto \frac{1}{\sigma_1^2}
\end{equation}
 Considering the probability density function of $[Z_1|\beta_1,\sigma_1^2,\theta_1] $ and the Bayes formula,
the posterior distribution of $[\beta_1|z_1,\sigma_1^2,\theta_1]$ is :
\begin{equation}
[\beta_1|z_1,\sigma_1^2,\theta_1] \sim \mathcal{N}_{p_1}\left(\lbrack F_1^T R_1(D_1)^{-1}F_1 \rbrack^{-1} \lbrack F_1^T R_1(D_1)^{-1}z_1\rbrack , \lbrack F_1^T \frac{R_1(D_1)^{-1}}{\sigma_1^2} F_1 \rbrack ^{-1 }\right)
\end{equation}
Then,  using the Bayes formula, we obtain that the posterior distribution of $[\sigma_1^2|z_1,\theta_1]$ is:
\begin{equation}\label{bayes36}
[\sigma_1^2|z_1,\theta_1] \sim \mathcal{IG}(\alpha_{\sigma_1^2|n_1}, \frac{Q_1}{2})
\end{equation}
where $ \mathcal{IG}(\alpha,Q)$ stands for the inverse gamma distribution with density function:
\begin{displaymath}
p_{\alpha,Q}(x) = \frac{Q^\alpha}{\Gamma(\alpha)}\frac{e^{-\frac{Q}{x}}}{x^{\alpha+1}}1_{x >0}
\end{displaymath}
and the parameters are given by:
\begin{equation}
\alpha_{\sigma_1^2|n_1} =
\frac{n_1-p_1}{2}
\qquad
\begin{array}{rl}
Q_1 & = (z_1 - F_1 \hat{\beta_1})^TR_1(D_1)^{-1}(z_1 - F_1 \hat{\beta_1}) \\
\end{array}
\end{equation}
with $
\hat{\beta_1} = E[\beta_1|z_1,\sigma_1^2,\theta_1] = \lbrack F_1^T R_1(D_1)^{-1}F_1 \rbrack^{-1} \lbrack F_1^T R_1(D_1)^{-1}z_1\rbrack
$.

Bayesian estimation of parameters with non-informative ``Jeffreys priors" \cite{J61} gives the same results as maximum likelihood estimation for the parameter $\beta_1$. For the parameter $\sigma_1^2$, the estimation given by $\frac{Q_1}{2\alpha_{\sigma_1^2|n_1}}$ is identical to the one obtained with the restricted maximum likelihood method. This method was introduced by Patterson and Thompson \cite{PT71} in order to reduce the bias of the maximum likelihood estimator. 

Secondly, let us consider the set of parameters $(\beta_2,\rho,\sigma_2^2,\theta_2)$. In order to have closed form formulas for the estimation of $(\beta_2,\rho)$, we   estimate them together. The idea to carry out a joint estimation is proposed for the first time in this paper and we believe it is important. Indeed, if the cheaper code is perfectly known, it can be considered as a regression function and so $\rho$ will be a regression parameter. In this case, it is clear that a separated estimation of $\beta_2$ and $\rho$ cannot be optimal.\\
Using similar Jeffrey prior distributions as in (\ref{JeffPrior01}) and the same methodology as for the estimation of $(\beta_1,\sigma_1^2)$, we find that:
\begin{equation}
[(\rho,\beta_2)|z_1,z_2,\sigma_2^2,\theta_2] \sim \mathcal{N}_{p_2+1}\left(\lbrack F^T R_2(D_2)^{-1}F\rbrack^{-1} \lbrack F^T R_2(D_2)^{-1}F\rbrack ,\lbrack F^T \frac{R_2(D_2)^{-1}}{\sigma_2^2}F\rbrack^{-1}\right)
\end{equation}
and:
\begin{equation}
[\sigma_2^2|z_2,z_1,\theta_2] \sim \mathcal{IG}(\alpha_{\sigma_2^2|n_2}, \frac{Q_2}{2})
\end{equation}
where:
\begin{equation}
\alpha_{\sigma_2^2|n_2} =
\frac{n_2-p_2-1}{2}
\qquad
\begin{array}{rl}
Q_2 &  = (z_2 - F \hat{\lambda})^TR_2(D_2)^{-1}(z_2 - F \hat{\lambda}) \\
\end{array}
\end{equation}
with $
\hat{\lambda} = E[(\rho,\beta_2)|z_1,z_2,\sigma_2^2,\theta_2] = \lbrack F^T R_2(D_2)^{-1}F\rbrack^{-1} \lbrack F^T R_2(D_2)^{-1}z_2\rbrack
$.
The design matrix $F$ is such that $F=[\rho z_1(D_2) \quad F_2] $. Furthermore, the estimation of $\sigma_2^2$ given by $\frac{Q_2}{2\alpha_{\sigma_2^2|n_2}}$ is  the same as the restricted maximum likelihood one. 

The hyper-parameters $\theta_1$ and $\theta_2$ are found by minimizing the opposite of the concentrated restricted log-likelihoods:
\begin{equation}\label{MLL01}
\mathrm{log}\left(\vert \mathrm{det}\left(R_1(D_1) \right) \vert \right)+ (n_1-p_1)\mathrm{log}(\hat{\sigma_1}^2)
\end{equation}
\begin{equation}\label{MLL02}
\mathrm{log}\left(\vert \mathrm{det}\left(R_2(D_2) \right) \vert \right)+ (n_2-p_2-1)\mathrm{log}(\hat{\sigma_2}^2)
\end{equation}
These minimizations problems must be numerically solved  with a global optimization method. We use an evolutionary method coupled with a BFGS algorithm. The drawback of the maximum likelihood estimation (see \cite{LC98}) is that, contrarily to Bayes estimation, we do not have any information about the variance of the estimator. Nevertheless, Bayes estimation of the hyper parameters $\theta_1$ and $\theta_2$ are prohibitive and as noted in \cite{S03} the choice of the prior distribution is non trivial. Therefore, in this paper,  we will always estimate these parameters  with a concentrated restricted likelihood method.

\section{Bayesian prediction for a code with 2 levels}\label{Baypred2levels}

The aim of a Bayesian prediction is to provide a predictive distribution for $Z_s(x)$ integrating the posterior distributions of the parameters and hence taking into account their uncertainty. 

A Bayesian prediction for a code with $s=2$ levels was suggested in  \cite{QW07}. Nevertheless, we propose here a new Bayesian approach with some significant differences. First, we assume that the adjustment coefficient is a regression function whereas  Qian and Wu \cite{QW07} model it with a Gaussian process. Secondly, we use different prior distributions for the parameter estimation. More specifically, according to the Bayesian estimation of parameters previously presented, we use a joint prior distribution for $(\beta_2,\rho)$ conditioned by  $\sigma_2^2$ whereas in \cite{QW07} they use separated prior distributions with $\rho$  not conditioned by $\sigma_2^2$. Then, we use a hierarchy between the different parameters. At the lowest level is the regressor parameter $\beta$. At the second level is the variance parameter $\sigma^2$ which controls the distribution of the parameter $\beta$. At the top level is the parameter $\theta$ which controls the distribution of the parameters at the bottom levels. It is common to use a hierarchical specification of models for Bayesian prediction as presented in \cite{R06}. This strategy will allow us to obtain explicit formulas for the joint distribution of the parameters and above all, to reduce dramatically  the cost of the numerical implementation of the complete Bayesian prediction.

We will also present the case in which we do not have any prior information about the parameters. As described in the previous section,  the hyper parameter $\theta$  is estimated by minimizing the opposite of the concentrated restricted log-likelihood and it is assumed to be fixed to this estimated value from now on.

\subsection{Prior distributions and Bayesian estimation of the parameters}

Many choices of priors can be made for the Bayesian modelling. Here we study the two following cases:
\begin{itemize}
\item[(I)] Priors for each parameter are informative.
\item[(II)] Priors for each parameter are non-informative.
\end{itemize}
For the non-informative case (II), we  use the improper distributions corresponding to the ``Jeffreys prior'' and then the posterior distributions are given in Section \ref{Bayesian_estimation_of_parameters}. Note that non-informative distributions are used when we do not have prior knowledge. For the informative case (I), we will consider the following prior distributions:
\begin{displaymath}
[\beta_1|\sigma_1^2] \sim \mathcal{N}_{p_1}(b_1,\sigma_1^2V_1), \qquad [(\rho,\beta_2)|z_1,\sigma_2^2] \sim \mathcal{N}_{1+p_2}\left(b_{\lambda} = \left(\begin{array}{c}b_\rho \\ b_2 \end{array} \right),\sigma_2^2V_{\lambda}=\sigma_2^2 \left(\begin{array}{cc}V_\rho & 0 \\ 0 & V_2 \\ \end{array} \right) \right)
\end{displaymath}
\begin{displaymath}
[\sigma_1^2] \sim \mathcal{IG}(\alpha_1,\gamma_1), \qquad [\sigma_2^2|z_1] \sim \mathcal{IG}(\alpha_2,\gamma_2)
\end{displaymath}
where $b_1 \in \mathbb{R}^{p_1}$, $b_{\lambda} \in \mathbb{R}^{1+p_2}$, $V_1$ is a $(p_1 \times p_1)$ diagonal matrix, $V_{\lambda}$ is a $((1+p_2) \times (1+p_2))$ diagonal matrix and $\alpha_1,\gamma_1,\alpha_2,\gamma_2 > 0$. The forms of the priors are chosen in order to be able to get closed form expressions for the posterior distributions. Note that there are enough free parameters in the priors to allow the user to prescribe their means and variances.
From the previous prior definitions, the posterior distributions of the parameters are:
\begin{equation}\label{eq8}
[\beta_1|z_1,\sigma_1^2] \sim \mathcal{N}_{p_1}(A_i^1 \nu_i^1, A_i^1) \qquad
[(\rho,\beta_2)|z_1,z_2,\sigma_2^2] \sim \mathcal{N}_{p_2+1}(A_i^{\lambda} \nu_i^{\lambda}, A_i^{\lambda})
\end{equation}
where:
\begin{equation}\label{eq8bis}
A_i^1 = \left\{
\begin{array}{lr} 
\lbrack F_1^T \frac{R_1^{-1}(D_1)}{\sigma_1^2}F_1+\frac{V_1^{-1}}{\sigma_1^2}\rbrack^{-1} & \quad i= (\mathrm{I}) \\ 
\lbrack F_1^T \frac{R_1^{-1}(D_1)}{\sigma_1^2} F_1 \rbrack^{-1} & \quad i=(\mathrm{II}) \\
\end{array}
\right.
\nu_i^1 = \left\{
\begin{array}{lr} 
\lbrack F_1^T \frac{R_1^{-1}(D_1)}{\sigma_1^2}z_1+\frac{V_1^{-1}}{\sigma_1^2} b_1 \rbrack & \quad i= (\mathrm{I}) \\
\lbrack F_1^T \frac{R_1^{-1}(D_1)}{\sigma_1^2}z_1\rbrack & \quad i= (\mathrm{II}) \\
\end{array}
\right.
\end{equation}
\begin{equation}\label{eq8ter}
A_i^{\lambda} = \left\{
\begin{array}{lr} 
\lbrack F^T \frac{R_2^{-1}(D_2)}{\sigma_2^2}F+\frac{V_{\lambda}^{-1}}{\sigma_2^2}\rbrack^{-1} & \quad i=(\mathrm{I}) \\ 
\lbrack F^T \frac{R_2^{-1}(D_2)}{\sigma_2^2} F \rbrack^{-1} & \quad i=(\mathrm{II}) \\
\end{array}
\right.
\nu_i^{\lambda} = \left\{
\begin{array}{lr} 
\lbrack F^T \frac{R_2^{-1}(D_2)}{\sigma_2^2}z_2+\frac{V_{\lambda}^{-1}}{\sigma_2^2} b_{\lambda} \rbrack & \quad i=(\mathrm{I}) \\
\lbrack F^T \frac{R_2^{-1}(D_2)}{\sigma_2^2}z_2\rbrack & \quad i= (\mathrm{II}) \\
\end{array}
\right.
\end{equation}
and $F=[\rho z_1(D_2) \quad F_2] $. Furthermore, we have:
\begin{equation}\label{eq9}
[\sigma_1^2|z_1] \sim \mathcal{IG}(\alpha_i^{\sigma_1^2|n_1}, \frac{Q_i^1}{2}), \qquad [\sigma_2^2|z_2,z_1] \sim \mathcal{IG}(\alpha_i^{\sigma_2^2|n_2}, \frac{Q_i^2}{2})
\end{equation}
where:
\begin{displaymath}
Q_i^1 = \left\{
\begin{array}{lr}
\gamma_1+(b_1-\hat{\beta_1})^T(V_1+[F_1^TR_1^{-1}(D_1)F_1]^{-1})^{-1}(b_1-\hat{\beta_1})+Q_2^1 & \quad i=(\mathrm{I}) \\
z_1 ^T[R_1^{-1}(D_1)-R_1^{-1}(D_1)F_1(F_1^TR_1^{-1}(D_1)F_1)^{-1}F_1^TR_1^{-1}(D_1)]z_1 & \quad i=(\mathrm{II}) \\
\end{array}
\right.
\end{displaymath}
\begin{displaymath}
Q_i^2 = \left\{
\begin{array}{lr}
\gamma_2+(b_{\lambda}-\hat{\lambda})^T(V_{\lambda}+[F^TR_2^{-1}(D_2)F]^{-1})^{-1}(b_{\lambda}-\hat{\lambda})+Q_2^2 & \quad i=(\mathrm{I}) \\
z_2 ^T[R_2^{-1}(D_2)-R_2^{-1}(D_2)F(F^TR_2^{-1}(D_2)F)^{-1}F^TR_2^{-1}(D_2)]z_2 & \quad i=(\mathrm{II}) \\
\end{array}
\right.
\end{displaymath}
\begin{displaymath}
\hat{\beta_1}=(F_1^TR_1^{-1}(D_1)F_1)^{-1}F_1^TR_1^{-1}(D_1)z_1 \qquad \hat{\lambda}=(F^TR_2^{-1}(D_2)F)^{-1}F^TR_2^{-1}(D_2)z_2
\end{displaymath}
\begin{displaymath}
\alpha_i^{\sigma_1^2|n_1} = \left\{
\begin{array}{lr}
\frac{n_1}{2}+\alpha_1 & \quad i=(\mathrm{I}) \\
\frac{n_1-p_1}{2} & \quad i= (\mathrm{II}) \\
\end{array}
\right.
\qquad
\alpha_i^{\sigma_2^2|n_2} = \left\{
\begin{array}{lr}
\frac{n_2}{2}+\alpha_2 & \quad i=(\mathrm{I}) \\
\frac{n_2-p_2-1}{2} & \quad i=(\mathrm{II}) \\
\end{array}
\right.
\end{displaymath}

Mixing of informative and non-informative priors are of course possible and easy to implement. As we will discuss in Subsection \ref{DiscCompExistMeth}  and see in the examples of Section  \ref{Toy_examples}, the use of informative priors has minor impact on the mean estimation but may have a strong impact  on variance estimation.

\subsection{Predictive distributions when $\beta_2, \rho,\sigma_1^2$ and $\sigma_2^2$ are known}

As a preliminary step towards the Bayesian prediction carried out in the next subsection, we give here  Bayesian prediction in the form of closed form expressions when the parameters $\beta_2$, $\rho$, $\sigma_1^2$ and $\sigma_2^2$ are known. The conditional distribution of $[Z_2(x)|Z=z,\beta_2,\rho,\sigma_1^2,\sigma_2^2] $ is given by:
\begin{equation}\label{eq10}
[Z_2(x)|\mathcal{Z}=z,\beta_2,\rho,\sigma_1^2,\sigma_2^2] \sim \normal{\mu_i(x)}{\sigma^2_i(x)}
\end{equation}
where:
\begin{displaymath}
\mu_i(x) = h'(x)^T \left(\begin{array}{c} A_i^1\nu_i^1\\ \beta_2 \\ \end{array} \right)
+ t(x)^TV^{-1}\left(z-H \left(\begin{array}{c} A_i^1\nu_i^1\\ \beta_2 \\ \end{array} \right) \right)
\end{displaymath}
\begin{displaymath}
\sigma^2_i(x) = s^2_{Z_2}(x) + k_1 A_i^1 k_1^T
\end{displaymath}
and $A_i^1$ and $\nu_i^1$ are defined by (\ref{eq8bis}). Note that the estimated variance is augmented by the term  $ k_1 A_i^1 k_1^T$ which quantifies the uncertainty due to the estimation of $\beta_1$.
$k_1$ is a $(1 \times p_1)$ vector composed of the $p_1$ first elements of the $(1\times p_1,1\times p_2)$ vector $k = (k_1 , k_2) = h'(x)^T - t(x)^T V^{-1}H$.
$H$ is given by (\ref{equationH}). The existence of closed form formulas is important as it will allow for a fast numerical implementation.

\subsection{Bayesian prediction}\label{Bayesian_prediction}

Before performing the Bayesian prediction we note that -  thanks to the explicit joint prior distribution for $\beta_2$ and $\rho$, the independence hypotheses and the hierarchical specification of the parameters -  conditioning on $\theta$, we have an explicit formula for the following joint density:

\begin{equation}
p(\beta_1,\beta_2,\rho,\sigma_1^2,\sigma_2^2|z_1,z_2)=
p(\beta_1|\sigma_1^2,z_1)
p(\beta_2,\rho|\sigma_2^2,z_1,z_2)
p(\sigma_1^2|z_1)
p(\sigma_2^2|z_1,z_2)
\end{equation}

This explicit joint density is an original result which contrasts with \cite{QW07} and which allows us to avoid prohibitive implementation for the Bayesian analysis.

First, we consider the predictive distribution with $\sigma_1^2$ and $\sigma_2^2$ known. Considering the conditional independence assumption between $(\delta(x))_{x\in Q}$ and $(Z_1(x))_{x\in Q}$, the probability density function of $[Z_2(x)| \mathcal{Z}=z, \sigma_1^2,\sigma_2^2]$ can be deduced from the following integral:
\begin{equation}\label{eq11}
p(z_2(x)| z_1, z_2,\sigma_1^2,\sigma_2^2) = \int_{\mathbb{R}^{1+p_2}}{ p(z_2(x)|z_1, z_2,\beta_2,\rho,\sigma_1^2,\sigma_2^2)p(\rho,\beta_2|z_1, z_2,\sigma_2^2) \, d\rho d\beta_2 }
\end{equation}
where $p(z_2(x)|z_1, z_2,\beta_2,\rho,\sigma_1^2,\sigma_2^2)$ is given by (\ref{eq03bis}). This integral has to be numerically evaluated. Since $[\rho,\beta_2|z_1,z_2,\sigma_2^2]$ has a known normal distribution given by (\ref{eq8}), we  here use a crude Monte-Carlo algorithm when the dimension of $\beta_2$ and $\rho$ is high, or a trapezoidal quadrature method when it is low.

Then, we infer from the parameters $\sigma_1^2$ and $\sigma_2^2$. Due to the independence between $(\delta(x))_{x\in Q}$ and $(Z_1(x))_{x\in Q}$, the probability density function of $[Z_2(x)| \mathcal{Z}=z]$ is:
\begin{equation}\label{eq12}
p(z_2(x)| z_1, z_2) = \int_{\mathbb{R}^{2}}{p(z_2(x)| z_1, z_2,\sigma_1^2,\sigma_2^2)p(\sigma_1^2| z_1)p(\sigma_2^2| z_1, z_2)
\, d\sigma_1^2 d\sigma_2^2}
\end{equation}
where $p(\sigma_1^2| z_1)$ and $p(\sigma_2^2| z_1, z_2)$ are given by (\ref{eq9}). This integral has also to be numerically evaluated. Since we have a double integration, a quadrature method will be efficient. We use here a trapezoidal numerical integration, defining the region of integration $[ \sigma_{1_{inf}}^{2} , \sigma_{1_{sup}}^{2} ] \times [ \sigma_{2_{inf}}^{2} , \sigma_{2_{sup}}^{2} ] $ from the equation (\ref{eq9}) and such that $p(\sigma_{1_{inf}}^{2} | z_1)$, $p(\sigma_{1_{sup}}^{2}| z_1)$ $p(\sigma_{2_{inf}}^{2}| z_1, z_2)$ and $p(\sigma_{2_{sup}}^{2}| z_1, z_2)$ are close to $0$. This region  essentially contains the support of the function. Furthermore, we create a non-uniform integration grid distributed with a geometric progression.

Finally $p(z_2(x)| z_1, z_2)$ is a predictive density function integrating the posterior distribution of parameters $(\beta_2,\rho, \beta_1,\sigma_1^2,\sigma_2^2)$. We hence have a predictive distribution taking into account the uncertainties due to the parameter estimations.

\subsection{Discussion about the numerical evaluations of the  integrals}\label{DiscCompExistMeth}

We saw in the previous section that we can obtain an analytical prediction when $\beta_2$, $\rho$, $\sigma_1^2$ and $\sigma_2^2$ are known. From this analytical formula, we can have a Bayesian prediction with only two nested  integrations. One of them can be approximated with a quadrature or a crude Monte Carlo method, which is not too expensive. The other is a double integration approximated with a quadrature method which is efficient and not expensive. Therefore, we do not use any Markov chain Monte Carlo method and we considerably reduce the time and the complexity of the method. This allows us to easily build an accurate Bayesian metamodel. Practically, we use 441 integration points to approximate (\ref{eq12}) and 1000 Monte-Carlo particles to approximate (\ref{eq11}). Therefore, we have  441000 call to the predictive density function  (\ref{eq10}). 

To avoid a prohibitive implementation, another approach has also been proposed in \cite{CG09}. They adopt a Bayes linear formulation which requires only the specification of the means, variances, and covariance. See \cite{GW07} for further details about the Bayes linear approach. The strength of this method is that its computational cost is low. Nonetheless, since it only focuses on posterior means and covariances, it does not provide the full posterior predictive distribution. 

Finally, we highlight the fact that our Bayesian procedure can be used to perform multi-fidelity analysis with more than 2 levels of code whereas the cost of the one presented in \cite{QW07} is too high to allow for such analysis. We illustrate in Section \ref{Exam2}  through an industrial case the great practical  importance of using more than 2 levels of code.

\section{Toy examples}\label{Toy_examples}

We will present in this section some co-kriging metamodels using one-dimensional functions inspired by the example presented in \cite{AF07}. For the following examples, we will use a non-Bayesian co-kriging model - \emph{i.e.} the one presented in \cite{KO00} - but with a Bayesian estimation of the parameters (see Section \ref{Bayesian_estimation_of_parameters}) and for the second example we will  use a Bayesian co-kriging. 

Furthermore, the correlation kernels are assumed to be:
\begin{displaymath}
r_t(x_i^{(k)} - x_j^{(l)} ; \theta_t) = \mathrm{exp}\left(- \frac{\| x_i^{(k)} - x_j^{(l)} \|^2}{\theta_t^2}\right) 
\end{displaymath}
where $
t,k,l=1,2 \quad 1\leq i\leq n_1 \quad 1\leq j \leq n_2$.

\begin{figure}[h]
\begin{center}
\includegraphics[width = 8 cm]{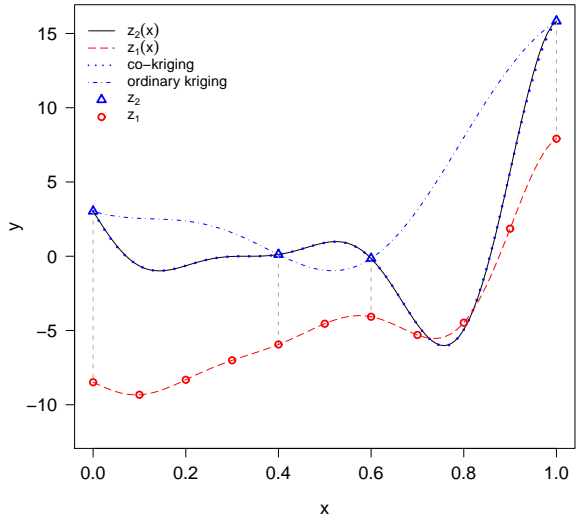}
\caption{Toy examples. The co-kriging metamodel is very close to the expensive output $z_2(.)$ and improves significantly the ordinary kriging metamodel using the small design $D_2$.}
\label{fig:academic01}
\end{center}
\end{figure}

\underline{Example 1.} The aim of this example is to emphasize the effectiveness of the presented Bayesian estimation of the parameters (see Section \ref{Bayesian_estimation_of_parameters}). We assume that the expensive code is given by $z_2(x) = (6x-2)^2 \mathrm{sin}(12x-4)$ and the cheaper code by $z_1(x) = 0.5z_2(x)+10(x-0.5)-5$. The experimental design set of the cheapest code is $D_1 = \{0,0.1,0.2,0.3,0.4,0.5,0.6,0.7,0.8,0.9,1 \}$ and the one of the expensive code is $D_2 = \{ 0,0.4,0.6,1 \}$. This example is identical to the one-dimensional demonstration presented in \cite{AF07}.
Figure \ref{fig:academic01} shows the functions $x \mapsto z_2(x)$ and $x \mapsto z_1(x)$, the training data for $z_2$ and $z_1$, the ordinary kriging using only the expensive data and the co-kriging using expensive and cheap data. 
To validate the model, the Root-mean-square errors (RMSE) and $Q_2 = 1 - \frac{\sum_{x \in T}{\left( m_{Z_2}(x)-z_2(x) \right)^2}}{\sum_{x \in T}{\left( m_{Z_2}(x)-\bar{z}_2 \right)^2}}$  are computed.

The test set $T$ is composed of a regular grid points sampled from $0$ to $1$ with a grid step equal to $0.01$ and $\bar{z}_2$ is the empirical mean evaluated in $T$. The estimated RMSE is $5.68 \times 10^{-2}$ and the coefficient $Q_2$ is $99.98 \%$, so we have a prediction error closed to 0. The Bayesian estimation of the parameters of co-kriging are given in Table \ref{tab1}. Furthermore, the estimations of the hyper-parameters $(\theta_1,\theta_2)$, calculated by maximizing the concentrated log-likelihoods (\ref{MLL01}) and (\ref{MLL02}), are $\hat{\theta}_1 = 0.25$ and $\hat{\theta}_2 = 0.80$. $D_1$ being a regular grid with a grid step equal to 0.1 and $D_2$ being composed of points sampled from 0 to 1, points of the experimental designs are hence strongly correlated which will imply a smooth surrogate model.
\begin{table}[H]
\begin{center}
\begin{tabular}{|c|c||c|c|}
\hline
Regression Coefficient & Estimation & Variance Coefficient & Estimation \\
\hline
$\rho$ & $2$ &$\sigma_1^2$ & $32.75$ \\
$\beta_2$ & $(20,-20)$ & $\sigma_2^2$ & $7.02\times10^{-30}$ \\
$\beta_1$ & $-3.49$ & & \\
\hline
\end{tabular}
\end{center}
\caption{A co-kriging example with one-variable functions. Bayesian estimation of parameters. }
\label{tab1}
\end{table}
We see that the Bayesian estimation of parameters is very effective since the estimations of parameters $\rho$ and $\beta_2$ are perfect. Nevertheless this example does not highlight the strength of the method since there is a relation between $z_2(x)_{x \in [0,1]}$ and $z_1(x)_{x \in [0,1]}$ which exactly corresponds to the equation (\ref{eq01}) with the error $\delta_2$ that can be written in terms of the regression functions $f_2$ exactly. Therefore, if the cheap code is well modelled, like in our case, the co-kriging is equivalent to a linear regression. Moreover, the very small value of $\sigma_2^2$ illustrates this.

\begin{figure}[h]
\begin{center}
\includegraphics[width = 8cm]{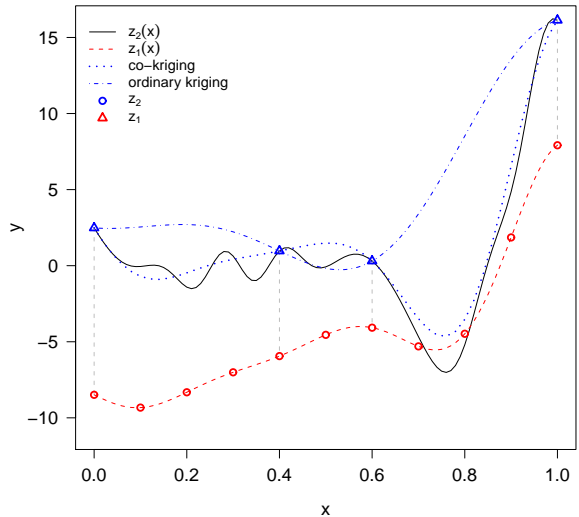}
\caption{Toy examples. The high frequency components of the expensive code are not predicted since they are not captured by the cheap code and the coarse grid used for the expensive code cannot detect them either. Nevertheless, the co-kriging improves the ordinary kriging metamodel since the cheap code allows us to predict the low frequencies of the expensive code accurately.}
\label{fig:academic02}
\end{center}
\end{figure}

\underline{Example 2.} This example illustrates a case where the non-Bayesian co-kriging underestimates the predictive variance whereas the Bayesian one adjusts it.  We assume that the expensive code is given by $z_2(x) = (6x-2)^2 \, \mathrm{sin} (12x-4) +\, \mathrm{sin}(10\,\mathrm{cos}(5x))$ and the cheaper code is given by $z_1(x) = 0.5( (6x-2)^2\,\mathrm{sin}(12x-4))+10(x-0.5)-5$. Through the term $\mathrm{sin}(10\,\mathrm{cos}(5x))$, the expensive code has high frequencies which are not captured by the cheap code and the error $\delta_2$ is not a simple linear combination of the regression functions $f_2$. Figure \ref{fig:academic02} shows the results of kriging and co-kriging for these two functions.
The estimated RMSE is $1.05$ and the coefficient $Q_2$ is $93.57\%$, we still have a good prediction. The Bayesian estimations of the parameters are given in Table \ref{tab2} and we have $\hat{\theta}_1 = 0.25$ and $\hat{\theta}_2 = 0.07$. The values of these parameters have been fixed according the following arguments. As the cheap code is the same as the one of the Example 1, we keep the same estimation for $\theta_1$. Then, we consider that there are not enough points to carry out a significant estimation of $\theta_2$. Therefore, we fix the value of $\hat{\theta}_2$ according to the high frequencies introduced by the term $\mathrm{sin}(10\,\mathrm{cos}(5x))$.
\begin{table}[H]
\begin{center}
\begin{tabular}{|c|c||c|c|}
\hline
Regression Coefficient & Estimation & Variance Coefficient & Estimation \\
\hline
$\rho$ & $1.86$ & $\sigma_1^2$ & $ 32.75.03$  \\
$\beta_2$ & $(18.39,-17.00)$ & $\sigma_2^2$ & $0.30$  \\
$\beta_1$ & $-3.49$ & & \\
\hline
\end{tabular}
\end{center}
\caption{A co-kriging example with one-dimensional functions. Bayesian estimation of parameters.}
\label{tab2}
\end{table}
Due to the additional term $\mathrm{sin}(10 \,\mathrm{cos}(5x))$, the estimation of the parameter $\rho$ is less effective than in the first example. This highlights the dependence between the estimation of $\rho$ and the mean of $\delta(x)_{x \in [0,1]}$.
\begin{figure}[h]
\begin{center}
\includegraphics[width = 8cm]{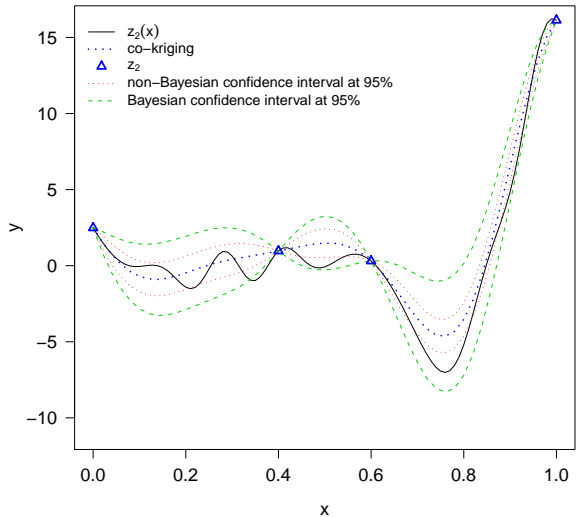}
\caption{A toy example without any prior information. The thick dotted line represents the prediction mean, the thin dotted lines represent the confidence interval at plus or minus twice the standard deviation in the non-Bayesian case and the dashed lines represent the same confidence interval in the Bayesian case.}
\label{fig:academic02bayes}
\end{center}
\end{figure}
Furthermore, Figure \ref{fig:academic02bayes} represents the confidence interval at plus or minus twice the standard deviation of the predictive distribution in the Bayesian and non-Bayesian case. We see that we underestimate the variance of the predictive distribution in the non-Bayesian case. This estimation is well adjusted in the Bayesian case.

\section{The case of $s$ levels of code}\label{CaseSLevels}

The aim of this Section is to perform a multi-level co-kriging with any number of codes. Let us consider $s$ levels of code. The generalization of the previous model is straightforward. Actually, if we note $\beta = (\beta_1^T,\dots,\beta_s^T)^T$, $\rho = (\rho_1,\dots,\rho_{s-1})$, $\sigma^2 = (\sigma_1^2,\dots,\sigma_s^2)$ and $\theta = (\theta_1,\dots,\theta_s)$, we have:
\begin{displaymath}
\forall x \in Q \qquad [Z_s(x)|\mathcal{Z}=z,\beta,\rho,\sigma^2,\theta] \sim \mathcal{N}\left(m_{Z_s}(x),s_{Z_s}^2(x)\right)
\end{displaymath}
where:
\begin{equation}\label{mZs}
m_{Z_s}(x) = h_s'(x)^T\beta + t_s(x)^TV_s^{-1}(z-H_s\beta)
\end{equation}
and:
\begin{equation}\label{s2Zs}
s_{Z_s}^2(x) = \sigma_{Z_s}^2 - t_s(x)^T V_s^{-1} t_s(x)
\end{equation}
Furthermore, let us denote by $R_t = R_t(D_t)$  the correlation matrix for $D_t$ and $\rho_s = 0$, $\forall s \leq 0$. The matrix $V_s$ has the form:
\begin{equation}\label{Vs}
V _s = \left(
\begin{array}{ccc}
V^{(1,1)} & \dots & V^{(1,s)}\\
\vdots & \ddots & \vdots \\
V^{(s,1)} & \dots & V^{(s,s)} \\
\end{array}
\right)
\end{equation}
The $s$ diagonal blocks of size $n_t \times n_t$ are defined by:
\begin{equation}
V ^{(t,t)} = \sigma_t^2 R_t(D_t) +\sigma_{t-1}^2 \rho_{t-1}^2 R_{t-1}(D_t) + \dots + \sigma_1^2 \left( \prod_{i=1}^{t-1} {\rho_i^2}\right)  R_1(D_t)
\end{equation}
and the off-diagonal blocks of size $n_t \times n_{t'}$ are given by:
\begin{equation}
V^{(t,t')} = \left(\prod_{i=t}^{t'-1} {\rho_i }\right)  V^{(t,t)}(D_t,D_{t'})  \qquad  1 \leq t < t' \leq s
\end{equation}
The vector $t_s(x)$ is such that $t_s(x) = (t_1^*(x,D_1)^T,\dots,t^*_s(x,D_{s})^T)^T$, where:
\begin{equation}
t^*_t(x,D_t)^T = \rho_{t-1} t^*_{t-1}(x,D_t)^T + \left(\prod_{i=t}^{s-1}{\rho_i}\right) \sigma_t^2R_t(x,D_t) \qquad
 1 < t \leq s 
\end{equation}
where $\left(\prod_{i=s}^{s-1}{\rho_i}\right) = 1$ and $t_1^*(x,D_1)^T = \left(\prod_{i=1}^{s-1}{\rho_i}\right)\sigma_1^2R_1(x,D_1)$. If we define:
\begin{displaymath}
F_k(D_l) = \left(
\begin{array}{c}
f_k^T(x_1^{(l)}) \\
\vdots \\
f_k^T(x_{n_l}^{(l)}) \\
\end{array}
\right)
\qquad
1 \leq k,l \leq s
\end{displaymath}
then the matrix $H_s$ can be written as:
\begin{equation}\label{Hs}
H_s = \left(
\begin{array}{llll}
F_1(D_1) & & &  \\
\rho_1 F_1(D_2) & F_2(D_2) & 0 & \\
\rho_1 \rho_2 F_1(D_3) & \rho_2 F_2(D_3) &  & \\
\vdots & \vdots &  \ddots & \\
 \left(\prod_{i=1}^{s-1} {\rho_{i}}\right)  F_1(D_s) & \left(\prod_{i=2}^{s-1} {\rho_{i} }\right)  F_2(D_s) & \dots & F_s(D_s)  \\
\end{array}
\right)
\end{equation}
$h_s'(x)^T$ and $\mathrm {var}(Z_s(x))=\sigma_{Z_s}^2$ are given by the equations (\ref{hsx}) and (\ref{s2sx}).

\subsection{Bayesian estimation of parameters for $s$ levels of code}\label{BayesEstimationslevels}

From the assumptions of conditional independence between $(\delta_t(x))_{x\in Q}$ and $(Z_{t-1}(x),\dots,Z_1(x))_{x\in Q}$, we can extend the Bayesian estimation of the parameters to the case of $s$ levels. Note that we do not assume the independence of $\beta_t$ and $\rho_{t-1}$. We can obtain a closed form expression for the estimation of $(\beta_t,\rho_{t-1})$. For all $t=2,\dots,s$, we have:
\begin{equation}
[(\rho_{t-1},\beta_t)|z_t,z_{t-1},\theta_t,\sigma_t^2] \sim \mathcal{N}\left( \left( H_t^T R_t^{-1}(D_t)H_t\right)^{-1} H_t^T R_t^{-1}(D_t)z_t, \sigma_t^2 \left( H_t^T R_t^{-1}(D_t)H_t\right)^{-1} \right)
\end{equation}
where $H_t = [\rho_{t-1}z_{t-1}(D_t)  \quad F_t(D_t)]$. Furthermore, if we note $
\hat{\lambda}_t = \mathbb{E}[(\rho_{t-1},\beta_t)|z_t,z_{t-1},\theta_t,\sigma_t^2] 
$,
then we have:
\begin{equation}
[\sigma_t^2|z_t,z_{t-1},\theta_t] \sim \mathcal{IG}(\alpha_t,\frac{Q_t}{2})
\end{equation}
where $
\alpha_t = \frac{n_t-p_t-1}{2} $ and $ Q_t = (z_t - H_t \hat{\lambda}_t)^T R_t^{-1}(D_t)(z_t - H_t \hat{\lambda}_t)
$.

The REML estimator of $\sigma_t^2$ is $ \hat{\sigma}_t^2 = \frac{Q_t}{2\alpha_t }$ and we can estimate $\theta_t$ by minimizing the expression:
\begin{equation}
\mathrm{log}(|\mathrm{det}(R_t(D_t))|) + (n_t-p_t-q_{t-1})\mathrm{log}(\hat{\sigma}_t^2 )
\end{equation}
The generalization of the Bayesian estimation previously presented is important since it shows that the parameters estimation for a $s$-levels co-kriging is equivalent to the one for $s$ independent krigings. 

\subsection{Reduction of computational complexity of inverting the covariance matrix $V_s$}\label{ImpRes}

$V_s$ is an $(\sum_{i=1}^{s}{n_i} \times \sum_{i=1}^{s}{n_i})$ matrix, its inverse can hence be difficult to process. We present in this Subsection a method to reduce the complexity of the processing of $V_s^{-1}$. By sorting the experimental design sets such that
$\forall t = 2,\dots,s $, $ D_{t-1} = (x_1^{(t-1)},\dots,x_{n_{t-1}-n_t}^{(t-1)}, x_1^{(t)},\dots,x_{n_t}^{(t)}) = (D_{t-1}\setminus D_t,D_t) 
$
it can be shown that $\forall t = 2,\dots,s$ the inverse of the matrix $V_s$ has the form:
\begin{equation}\label{big_inverse}
V_s^{-1} = \left(
\begin{array}{cc}
V_{s-1}^{-1} + \left(
\begin{array}{cc}
0 & 0 \\
0 &\rho_{s-1}^2 \frac{R_{s}^{-1}}{\sigma_{s}^2} \\
\end{array}
\right) &
- \left(
\begin{array}{c}
0 \\
\rho_{s-1} \frac{R_{s}^{-1}}{\sigma_{s}^2} \\
\end{array}
\right) 
\\
- \left(
\begin{array}{cc}
0 &
\rho_{s-1} \frac{R_{s}^{-1}}{\sigma_{s}^2} \\
\end{array}
\right) & 
\frac{R_{s}^{-1}}{\sigma_{s}^2} \\
\end{array}
\right)
\qquad
V_1^{-1} = \frac{R_{1}^{-1}}{\sigma_{1}^2}
\end{equation}
with $V_{s-1}^{-1}$ an $(\sum_{i=1}^{s-1}{n_i} \times \sum_{i=1}^{s-1}{n_i})$ matrix and $R_{s}^{-1}$ an $(n_s \times n_s)$ matrix. This is a very important result since it shows that we can deduce $V_s^{-1}$ from $R_{t}^{-1}$, $t=1,\dots,s$. Therefore, the complexity of the processing of $V_s^{-1}$ is $\mathcal{O}(\sum_{i=1}^{s}{n_i^3})$ instead of $\mathcal{O}((\sum_{i=1}^{s}{n_i})^3)$. Furthermore, from the equation (\ref{big_inverse}) and the Bayesian estimation of parameters presented in Section \ref{BayesEstimationslevels}, we have shown here that building a $s$-level co-kriging is equivalent to build $s$ independent krigings. 
We emphasize that, for practical applications, the form  (\ref{big_inverse}) for the inverse of $V_s$ allows us to perform fine matrix regularization in the case of ill-conditioned problems. Indeed, $V_s$ is invertible if and only if the matrices $R_t$, $t=1,\dots,s$ are invertible. Therefore, if the problem is ill-conditioned, we just have to regularize the matrices $R_t$ which are ill-conditioned too.
Then, since $ (t_1^*(x,D_1)^T,\dots,t^*_{s-1}(x,D_{s-1})^T) = \rho_{s-1} t_{s-1}^T(x) $ it can also be shown that in the equation (\ref{mZs}):
\begin{equation}
t_s(x)^T V_s^{-1} = \left(
\begin{array}{c@{,}c}
\rho_{s-1}t_{s-1}^T(x) V_{s-1}^{-1} -[0_{1\times ( \sum_{i=1}^{s-1}{n_{i}}-n_s)} , \rho_{s-1} R_s(\{x\},D_s)R_s^{-1} ] &
R_s(\{x\},D_s)R_s^{-1}
\end{array}
\right)
\end{equation}
Therefore, $t_s(x)^T V_s^{-1}$ is independent of $\sigma_s^2$. Since $t_1(x)^T V_1^{-1} = R_1(\{x\},D_1)R_1^{-1}$ does not depend on $\sigma_1^2$, by induction, $t_s(x)^T V_s^{-1}$ is independent of $\sigma_i^2$ for all $1\leq i \leq s$. We have just shown here that the co-kriging mean does not depend on the variance coefficients.

\subsection{Numerical test on the reduction of computational complexity}

In the previous section, we have  presented a reduction of complexity for the co-kriging model by expressing the inverse of the matrix $V_s$ with the inverses of the matrices $R_t$, $t=1,\dots,s$. We present here a numerical test to highlight the gain of CPU time obtained with this method. We focus on the case of 2 levels of code and we consider the Gaussian kernel for the 2 levels:

\begin{displaymath}
r(x-x';\theta)=\mathrm{exp}\left( - \frac{||x-x'||^2}{\theta^2}\right)
\end{displaymath}

\begin{figure}[h]
\begin{center}
\includegraphics[width =8cm]{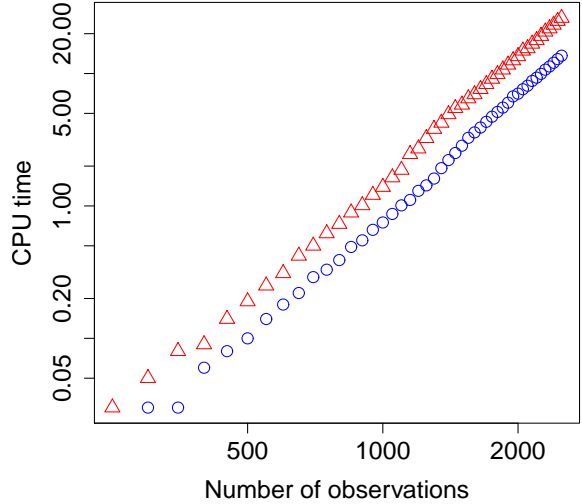}
\caption{CPU time comparison between 2 level co-kriging models. The triangles represent the CPU time for the crude co-kriging model and the circles represent the CPU time for the co-kriging model with the complexity reduction. The gain of CPU time  with the reduction complexity is approximately a factor equal to 1.93.}
\label{CPUtime}
\end{center}
\end{figure}

The experimental design set for the cheap code is a regular grid composed of $n_1$ points between 0 and 1 and the experimental design set for the expensive code are the $n_2$ first points of this grid. We consider the relation $n_1 = 4n_2$ with $n_2 = 50,60,\dots, 500$ and the parameter $\theta = \frac{5}{n_2}$ (the parameter $\theta$ is controlled by $n_2$ in order to avoid ill-conditioned covariance matrices). The total number of observations is hence $n = n_1 +n_2$. Figure \ref{CPUtime} compares the CPU time needed to build a co-kriging model with or without reduction complexity.

First, the slope of the two CPU times is close to 3  (the least-squares estimation value is  3.03). The complexity of a matrix inversion being  $\mathcal{O}(n^3)$, with $n$ the size of the matrix, the estimation of the slope highlights the fact that it is the matrix inversion which leads the CPU time. Then, Figure \ref{CPUtime} emphasizes that the reduction of complexity is worthwhile. Indeed, we see that  the ratio between the two CPU time is approximately a constant  equal to 1.93. We are hence close to the theoretical ratio   equal to $(n_1+n_2)^3/(n_1^3+n_2^3) \approx 1.92$ which is obtained when we consider that the CPU time is essentially due to the matrix inversion.

\subsection{Toy example on the complexity reduction}

A 3-level co-kriging metamodel is presented in this section to illustrate the gain of CPU which can be obtained with the presented reduction of complexity. We focus on the inversion of the co-kriging matrix $V_s$  by comparing the CPU time needed with a direct inversion or by using the formula (\ref{big_inverse}).
We assume that the 3 levels of code are given by the followings three dimensional functions:
\begin{eqnarray}
z_1(x) & = & \mathrm{sin}(x_1)\\
z_2(x) & = &z_1(x)+a\mathrm{sin}(x_2)^2 \\
z_3(x) & = & z_2(x)+b x_3^4\mathrm{sin}(x_1)
\end{eqnarray}
with $x = (x_1,x_2,x_3) \in [-\pi,\pi]^3$, $a = 7$ and $b = 1/10$. We note that the complex function $z_3(x)$ corresponds to the Ishigami function which is very popular in the field of sensitivity analysis \cite{Sal00}.
We consider $n_3 = 50$ observations for the most accurate code $z_3(x)$, $n_2 = 200$ for the intermediate code and $z_1 = 400$ for the less accurate code. All experimental design sets are randomly sampled from the uniform distribution. As presented in Section \ref{building} we consider nested experimental designs $\forall t=2,\dots,s \quad D_t \subseteq D_{t-1}$.

We use a tensorised Mat\'ern-$\frac{5}{2}$ kernel for the three correlation functions:
\begin{equation}
r_t(x,x';\theta_t) = \prod_{i=1}^{d}r_{1D}(x_i,x'_i;\theta_{t,i})
\end{equation}
with
$
r_{1D}(t,t';\theta) = \left(1+ \sqrt{5}\frac{|t-t'|}{\theta} + \frac{5}{3}\frac{(t-t')^2}{\theta^2}\right)\mathrm{exp}\left(-\sqrt{5}\frac{|t-t'|}{\theta} \right)
$, $t,t' \in \mathbb{R}$.

The estimations of the hyper-parameters $\theta_t$ are presented in Table \ref{hypparamishi}.
\begin{table}[H]
\begin{center}
\begin{tabular}{|c|c|}
\hline
\textbf{Parameter} & \textbf{Estimation} \\
\hline
$\hat{\theta}_1$ &  $(\begin{array}{ccc} 0.61 &  1.99 &  2.04  \end{array})$ \\
\hline
$\hat{\theta}_2$ & $(\begin{array}{ccc} 1.98 &  0.26 &  2.48 \end{array})$ \\
\hline
$\hat{\theta}_3$ & $(\begin{array}{ccc} 0.23 &  0.89 &  0.21 \end{array})$  \\
\hline
\end{tabular}
\end{center}
\caption{Toy example on the complexity reduction. Estimation of the hyper-parameters (correlation lengths) for the 3-level co-kriging.}
\label{hypparamishi}
\end{table}
The hyper-parameter estimates show us that $z_1(x)$ is very smooth in the directions $x_2$ and $x_3$ reflecting the fact that it depends only on the first direction $x_1$. Similarly, the bias between $z_2(x)$ and $z_1(x)$ only depending  on the second direction $x_2$, it is rougher on this direction and very smooth in the others. Finally, the bias between $z_3(x)$ and $z_2(x)$ is rougher in the  direction $x_3$ than in  the directions $x_1$ and $x_2$. This is due to the important impact of $x_3$ on the third level.

The estimation of the variance, scale and regression parameters are given in Table \ref{paramishi}.
\begin{table}[H]
\begin{center}
\begin{tabular}{|c|c|}
\hline
Parameter & Estimation \\
\hline
$\beta_1$ &0.00 \\
$\left(\begin{array}{c} \rho_1 \\ \beta_2 \\ \end{array} \right)$ & $ \left( \begin{array}{c} 0.99 \\ 2.44 \\ \end{array} \right)$ \\
$\left(\begin{array}{c} \rho_2 \\ \beta_3 \\ \end{array} \right)$ & $ \left( \begin{array}{c} 0.95 \\ 0.64 \\ \end{array} \right)$ \\
\hline
$\sigma_1^2$ & 0.09  \\
$\sigma_2^2$ & 1.66  \\
$\sigma_3^2$ & 6.25 \\
\hline
\end{tabular}
\end{center}
\caption{Toy example on the complexity reduction. Estimation of the variance, scale and regression parameters for the 3-level co-kriging.}
\label{paramishi}
\end{table}
Table \ref{paramishi} shows the efficiency of the suggested method for the parameter estimations since it provides very accurate estimations of $\rho_1$ and $\rho_2$. 

To evaluate the accuracy of the co-kriging model, we use a test set of 30,000 points uniformly sampled from the uniform distribution. Then,  we compute the coefficient $Q_2$ with the co-kriging predictions and the responses of $z_3(x)$ on this set. We obtain for the co-kriging model $Q_2 = 83.21\%$, we hence have a middling accuracy despite the large number of observations used. Nonetheless, we have a significant improvement relatively to the kriging model since with the same kernel and the same experimental design set $D_3$ we obtain   $Q_2 = 47.97\%$ which is  a very poor accuracy. The hyper-parameter estimation of the kriging model is $\theta = (0.79,0.14,0.29)$, the variance one is $\sigma^2 = 13.66$ and the trend coefficient one is $\beta = 3.89$.

Let us now compare the difference of CPU time between the co-kriging building with a crude inversion of the covariance matrix $V_s$ and the one with an inversion using  the formula presented in Subsection \ref{ImpRes}. The CPU time necessary without the reduction complexity is CPU$_\mathrm{crude}=0.47$ whereas the one necessary with the complexity reduction is CPU$_\mathrm{light}=0.14$. We hence find that the CPU time ratio  between the two methods approximately equals $3.36$. This is not far from the theoretical ratio which equals $650^3/(400^3+200^3+50^3) \approx 3.80$. We note that the complexity reduction could be  of important practical interest. For example, without it the computational cost of a leave-one-out cross validation procedure will be much more important (the ratio will still be around 3 in our example). The complexity of this procedure being $\mathcal{O}(n^4)$, the gain of CPU time will be substantially.

\section{Example : Fluidized-Bed Process}\label{Exam2}
This example illustrates the comparison between 2-level and 3-level co-kriging. A 3-level co-kriging method is applied  to a physical experiment modelled by a computer code. The experiment, which is the measurement of the temperature of the steady-state thermodynamic operation point for a fluidized-bed process, was presented by
\cite{Dew99}, who developed a computer model named ``Topsim'' to calculate the measured temperature. The code, developed for a Glatt GPCG-1, fluidized-bed unit in the top-spray configuration, can be run at 3 levels of complexity. We hence have 4 available responses:
\begin{enumerate}
\item $T_{exp}$: the experimental response.
\item $T_3$: the most accurate code modelling the experiment.
\item $T_2$: a simplified version of $T_3$.
\item $T_1$: the lowest accurate code modelling the experiment.
\end{enumerate}
The differences between $T_1$, $T_2$ and $T_3$ are discussed by Dewettinck et al. (1999). The aim of this study is to predict the experimental response $T_{exp}$ given the two levels of code $T_3$ and $T_2$. We only focus on a 3-level co-kriging using $T_3$ and $T_2$ to predict $T_{exp}$ since 28 responses available for each level is not enough  to build a nested experimental design relevant for a  4-level co-kriging. The experimental design set and the responses $T_1$, $T_2$, $T_3$ and $T_{exp}$ are given by \cite{QW07} who have presented a 2-level co-kriging using $T_{exp}$ and $T_2$. Furthermore, the responses are parameterized by a 6-dimensional input vector presented by Dewettinck et al. (1999). 
\subsection{Building the 3-level co-kriging}
To build the 3-level co-kriging, we use 10 measures of $T_{exp}$ (measures 1, 3, 8, 10, 12, 14, 18, 19, 20, 27 in Table 4 in \cite{QW07}), 20 simulations of $T_3$ (runs 1, 2, 3, 5, 6, 7, 8, 9, 10, 11, 12, 13, 14, 16, 18, 19, 20, 22, 24, 27) and the 28 simulations of $T_2$ and the input vector is scaled between 0 and 1. The last 18 measures of $T_{exp}$ are used for validation. The  design sets are nested such that $D_{t-1} = (D_{t-1} \setminus D_t , D_t)$ for $t=2,3$ and we use a Matern$\frac{5}{2}$ kernel for the three covariance functions. The estimations of the hyper-parameters which represent correlation lengths of the three covariance kernels are given in Table \ref{tab2lt}.
\begin{table}[H]
\begin{center}
\begin{tabular}{|c|c|c|c|c|c|c|}
\hline
$\hat{\theta}_1$ &1.790 &3.988& 1.218 &1.790& 3.595& 0.722 \\
\hline
$\hat{\theta}_2$ &1.810 &1.842& 2.008 &1.036& 0.001& 0.345 \\
\hline
$\hat{\theta}_3$ &0.890 &0.721& 2.008 &2.952&1.790& 0.241 \\
\hline
\end{tabular}
\end{center}
\caption{Example: fluidized-bed process. Estimation of the hyper-parameters (correlation lengths) for the 3-level co-kriging.}
\label{tab2lt}
\end{table}
The estimations of hyper-parameters in Table \ref{tab2lt} show us that the surrogate model will be very smooth in the first four directions. For the fifth direction the Gaussian processes modelling the cheap code $T_2$ and the bias between $T_{exp}$ and $T_3$ are very smooth and the one modelling the bias between $T_3$ and $T_2$ is close to a regression. Finally, the model is sharper in the sixth direction in particular for the two biases where correlation lengths are around 0.3.

Furthermore,  Table \ref{tab2lbs} gives the estimation of the variance and regression parameters (see section \ref{BayesEstimationslevels}).
\begin{table}[H]
\begin{center}
\begin{tabular}{|c|c|c|}
\hline
Regression coefficient & Posterior mean & $\frac{\mathrm{Posterior\ Covariance}}{\sigma_t^2}$ \\
\hline
$\beta_1$ & 47.02 & 0.134 \\
$\left(\begin{array}{c} \beta_{\rho_1} \\ \beta_2 \\ \end{array} \right)$ & $ \left( \begin{array}{c} 0.97 \\ -0.17 \\ \end{array} \right)$ & $ \left( \begin{array}{cc} 0.001 & -0.034 \\-0.034 &1.610 \\ \end{array} \right) $\\
$\left(\begin{array}{c} \beta_{\rho_2} \\ \beta_3 \\ \end{array} \right)$ & $ \left( \begin{array}{c} 0.95 \\ 1.93 \\ \end{array} \right)$ & $ \left( \begin{array}{cc} 0.003 &-0.121 \\ -0.121 & 5.188 \\ \end{array} \right) $ \\
\hline
Variance coefficient & $Q_t$ & $\alpha_t$ \\
\hline
$\sigma_1^2$ & 1032 & 13.5 \\
$\sigma_2^2$ & 5.30 & 9 \\
$\sigma_3^2$ & 8.39 & 4 \\
\hline
\end{tabular}
\end{center}
\caption{Example: fluidized-bed process. Bayesian estimation of the variance and regression parameters for the 3-level co-kriging.}
\label{tab2lbs}
\end{table}

  Table \ref{tab2lbs} shows that the responses have approximately the same scale since the adjustment coefficients are close to 1. Furthermore, we see an important bias between $T_3$ and $T_2$ with $\beta_3 = 1.93$. Finally,  the variance coefficients for the biases indicate that they are possibly much simpler to model than the cheap code $T_2$ as their estimations are smaller.

\subsection{3-level co-kriging prediction: predictions when code output is available}\label{build3lev}

The aim of this Section is to show that co-kriging can improve significantly the accuracy of the surrogate model at points where at least one level of responses is available.

\begin{figure}[h]
\begin{center}
\includegraphics[width =8cm]{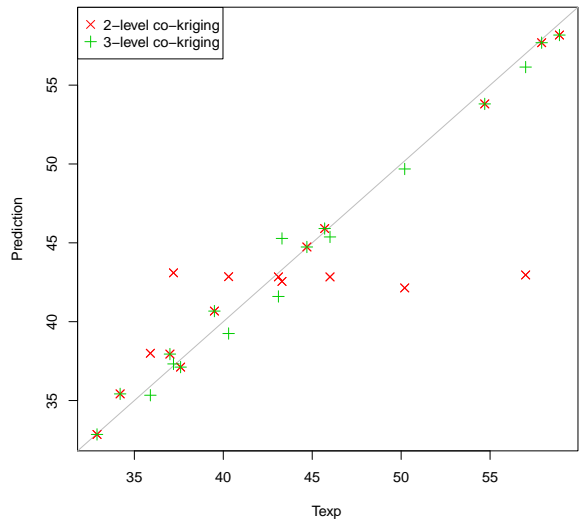}
\caption{Predictions of the 2-level and the 3-level co-krigings for the fluidized-bed process. The 3-level co-kriging improves significantly the predictions of the 2-level one.}
\label{fig:3lpred}
\end{center}
\end{figure}

\begin{figure}[h]
\begin{center}
\includegraphics[width =8cm]{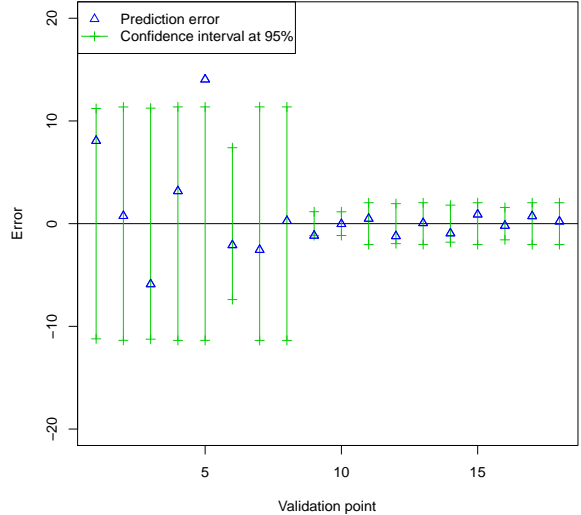}
\caption{Prediction errors of the 2-level co-kriging and confidence intervals at plus or minus twice the standard deviation. We see a significant difference between the accuracy of the predictions means and their confidence intervals for the point where $T_3$ is unknown (the 8 first validation points) and for the ones where it is known (the last 10 validation points).}
\label{fig:3lconf}
\end{center}
\end{figure}

The predictions of the 3-level co-kriging are here  presented  and compared with the predictions obtained with a 2-level co-kriging using only the 10 responses of $T_{exp}$ and the 20 responses of $T_3$. The predictions for the 2-level and the 3-level co-krigings vs. the real values (i.e., the measured temperature $T_{exp}$) are shown in Figure \ref{fig:3lpred}.
The 3-level co-kriging gives us the same prediction means as the 2-level co-kriging at the 10 points (points 2, 5, 6, 7, 9, 11, 13, 16, 22, 24) where $T_3$ is known. These overlapped points mean that $T_2$ does not influence the surrogate model at these points. This follows from the Markov property introduced in Section \ref{building}, which implies that the prediction of $T_{exp}$ is entirely determined by $T_3$ at these points. We also note that, in general, the 2-level co-kriging predictions - at points where  $T_3$ is unknown - are not accurate and  the 3-level co-kriging improves significantly the prediction means compared to the 2-level co-kriging. Table \ref{perf3l} compares the 2-level co-kriging with the 3-level co-kriging and summarizes some results about the quality of the predictions on the 18 validation points. Nonetheless, it is important to notice that, in the 3-level case, the output of the cheapest code $T_2$ is known at the 18  test points. This means that the results of this subsection show that the 3-level co-kriging prediction is more accurate than the 2-level co-kriging prediction at a point where the cheapest response $T_2$ is available. In the next subsection we show that the 3-level co-kriging prediction is more accurate than the 2-level one at a point where no response is available.
\begin{table}[H]
\begin{center}
\begin{tabular}{cccc}
\hline
& $Q_2$ & RMSE & MaxAE \\
2-level co-kriging & 61.23 \% & 4.24 & 14.04 \\
3-level co-kriging & 98.71 \% & 0.89 & 1.98 \\
\hline
\hline
& Average Std. dev. & Median Std. dev. & Maximal Std. dev \\
2-level co-kriging & 2.90 & 1.02 & 5.68  \\
3-level co-kriging & 0.90 & 1.02 &  1.04 \\
\hline
\end{tabular}
\end{center}
\caption{Example: fluidized-bed process. Comparison between 2-level co-kriging and 3-level co-kriging. Predictions are better in the 3-level case and the prediction variance seems well-evaluated since the RMSE and the average standard deviation are close.}
\label{perf3l}
\end{table}
Figure \ref{fig:3lconf} shows the prediction errors of the 2-level co-kriging and the confidence interval at plus or minus twice the prediction standard deviation. The last 10 prediction errors and their confidence intervals are the same as those of the 3-level case since it corresponds to the points where $T_3$ is known.
We see in Figure \ref{fig:3lconf} that the confidence intervals are well predicted. Furthermore, we see a significant difference between the accuracy of the prediction means and their confidence intervals for the point where $T_3$ is unknown (the 8 first validation points) and for the ones where it is known (the last 10 validation points).

\subsection{3-level co-kriging prediction: predictions when code output is not available}

In this subsection, we show that a multi-level co-kriging can significantly improve the prediction of a surrogate model at points where no response is available.

\begin{figure}[h]
\begin{center}
\includegraphics[width = 8cm]{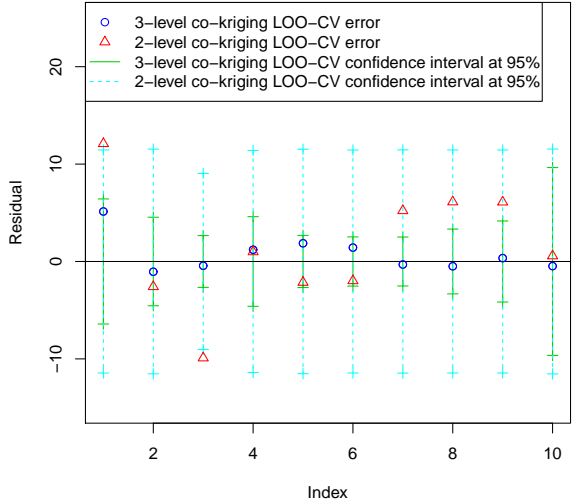}
\caption{Leave-One-Out Cross-Validation predictive errors and variances of the 2-level and 3-level co-kriging. We see that the confidence intervals are accurate and the precision of the 3-level co-kriging is significantly better than the one of the 2-level co-kriging.}
\label{fig:3l2lLOO}
\end{center}
\end{figure}

We have seen in Section \ref{build3lev} that the 3-level co-kriging improves significantly the 2-level co-kriging at points where $T_3$ is unknown and $T_2$ has been sampled. Nevertheless,  to have a fair comparison between these two co-kriging models, we  compare their accuracy by applying a Leave-One-Out Cross-Validation (LOO-CV) procedure at the 10 points where  $T_{exp}$ is known. This means that we perform for each of these 10 points the following procedure: 
\begin{enumerate}
\item The experimental and the two code outputs corresponding to the point are removed from the data set.
\item The 2-level co-kriging method and the 3-level co-kriging method are applied using the truncated data set in order to give a confidence interval for the experimental output at the point.
\end{enumerate}
Figure \ref{fig:3l2lLOO} shows the result of the LOO-CV procedure for the 2-level and 3-level co-kriging. We see that the 3-level co-kriging is more accurate than the 2-level one. Indeed, the LOO-CV RMSE for the 2-level co-kriging is equal to  1.88 whereas it is equal to 1.09 for the 3-level co-kriging. This shows that the 3-level co-kriging provides better predictions also at points where no response is available. This highlights the strength of the proposed method and shows that a co-kriging method with more than 2 levels of code can be worthwhile.

\section{Conclusion}

We have  presented a method for building kriging models using a hierarchy of codes with different levels of accuracy. This method allows us to improve a surrogate model built on a complex code using information from a cheap one. It is particularly useful when the complex code is very expensive. We see in our literature review that the first multi-level metamodel originally suggested is a first-order auto-regressive model built with Gaussian processes. The AR(1) relation between two levels of code is natural and the building of the model is straightforward. Nevertheless, we have highlighted some key issues which makes it difficult to use this model in practical ways.

First, important parameters of the model, which are the adjustment coefficients between two successive levels of codes, were numerically estimated. We propose here an analytical estimation of these parameters with a Bayesian method. This method allows us to have information about the uncertainties of the estimations and above all, to easily use the AR(1) model and its generalization to the case of non-spatial stationarity. Furthermore, a strength of the proposed method is that it even works for a code with more than 2 levels since its implementation is such that the estimations of the parameters of a $s$-level co-kriging is equivalent to  the ones of $s$ independent krigings. It is important to highlight that this method is based on a joint estimation between the adjustment coefficient and the mean of the Gaussian process modelling the difference between two successive levels of code.

Second, we have seen that the variance of the predictive distribution of the AR(1) model could be underestimated. A natural approach to improve this estimation is a Bayesian modelling. We propose here a Bayesian co-kriging for 2 levels of code and to avoid computationally expensive implementation, we suggest another model than the one  presented. This new model is based on a hierarchical specification of the parameters of the model. This allows us to have a Bayesian model including only two nested integrations without Markov chain Monte Carlo procedure.

Finally, for a non-Bayesian $s$-level  co-kriging, we have proved that building a $s$-level co-kriging is equivalent to build $s$ independent krigings. This result is very important since it solves one of the most important key issues of the co-kriging which is the inversion of the covariance matrix. A 3-level co-kriging example has been provided to show the efficiency of the presented method.

\section{Acknowledgements}

The author thanks his supervisor Josselin Garnier for his valuable guidance. He also thanks Claire Cannamela for her advice and constructive suggestions.

\appendix

\section{The case of $\rho$ depending on $x$}\label{appendix1}

\subsection{Building a model with $s$ levels of code}

 Let us consider $s$ levels of code, if we note $\beta = (\beta_1^T,\dots,\beta_s^T)^T$, $\beta_\rho = (\beta_{\rho_1}^T,\dots,\beta_{\rho_{s-1}}^T)^T$, $\sigma^2 = (\sigma_1^2,\dots,\sigma_s^2)$ and $\theta = (\theta_1,\dots,\theta_s)$, we have $[Z_s(x)|\mathcal{Z}=z,\beta,\beta_\rho,\sigma^2,\theta] \sim \mathcal{N}\left(m_{Z_s}(x),s_{Z_s}^2(x)\right)
$
where $m_{Z_s}(x)$ and  $s_{Z_s}^2(x)$ are defined in equations (\ref{mZs}) and (\ref{s2Zs}).
Let us define the  notation $\bigodot_{i=k}^l A_i = A_k \odot \dots \odot A_l$
where $\odot $ represents the matrix element-by-element product. Furthermore, let us denote by $\rho_t = \rho_t(D_t)$ the vector containing the values of $\rho_t(x)$, $x\in D_t$. The $s$ diagonal blocks of $V_s$  (\ref{Vs}) of size $n_t \times n_t$ are defined by:
\begin{displaymath}
V ^{(t,t)} = \sigma_t^2 R_t(D_t) +\sigma_{t-1}^2 \left( \rho_{t-1}(D_t) \rho_{t-1}^T(D_t)\right) \odot R_{t-1}(D_t) + \dots + \sigma_1^2 \left( \bigodot_{i=1}^{t-1} {\rho_i (D_t) \rho_i^T(D_t)}\right) \odot R_1(D_t)
\end{displaymath}
and the off-diagonal blocks of size $n_t \times n_{t'}$ are given by:
\begin{displaymath}
V^{(t,t')} = \left(\mathbf{1}_{n_{t}}\left(\bigodot_{i=t}^{t'-1} {\rho_i }(D_{t'}) \right)^T\right) \odot V^{(t,t)}(D_t,D_{t'})  \qquad 1 \leq t < t' \leq s
\end{displaymath}
The vector $t_s(x)$ in equations (\ref{mZs}) and (\ref{s2Zs})  is such that $t_s(x) = (t_1^*(x,D_1)^T,\dots,t^*_s(x,D_{s})^T)^T$, where:
\begin{displaymath}
t^*_t(x,D_t)^T = \rho_{t-1}^T(D_t) \odot t^*_{t-1}(x,D_t)^T + \left(\prod_{i=t}^{s-1}{\rho_i(x)}\right) \sigma_t^2R_t(x,D_t)
\end{displaymath}
where $ 1 < t \leq s$, $\left(\prod_{i=s}^{s-1}{\rho_i(x)}\right) = 1$ and $t_1^*(x,D_1)^T = \left(\prod_{i=1}^{s-1}{\rho_i(x)}\right)\sigma_1^2R_1(x,D_1)$. Furthermore, the matrix $H_s$ in equations \ref{Hs} can be written as:
\begin{displaymath}
H_s = \left(
\begin{array}{lllll}
&\vdots & \ddots  &  & \\
\left( \left(\bigodot_{i=1}^{j-1} {\rho_{i}(D_j) }\right)\mathbf{1}_{p_{1}}^T\right)\odot F_1(D_j) & \left(\left(\bigodot_{i=2}^{j-1} {\rho_{i}(D_j) }\mathbf{1}_{p_{2}}^T\right) \right) \odot F_2(D_j) & \dots & F_j(D_j) & 0  \\
& \vdots & &  & \ddots \\
\end{array}
\right)
\end{displaymath}

\subsection{Bayesian estimation of parameters for $s$ levels of code}\label{BayesEstimationslevelsA}

We can extend the Bayesian estimation of the parameters to the case of $\rho$ depending on $x$. Note that we do not assume the independence of $\beta_t$ and $\beta_{\rho_{t-1}}$. We have:
\begin{displaymath}
[(\beta_{\rho_{t-1}},\beta_t)|z_t,z_{t-1},\theta_t,\sigma_t^2] \sim \mathcal{N}\left( \left( H_t^T R_t^{-1}(D_t)H_t\right)^{-1} H_t^T R_t^{-1}(D_t)z_t, \sigma_t^2 \left( H_t^T R_t^{-1}(D_t)H_t\right)^{-1} \right)
\end{displaymath}
where $H_t = [F_{\rho_{t-1}}(D_t)\odot (z_{t-1}(D_t)\mathbf{1}_{q_{t-1}}^T)  \quad F_t(D_t)]$. Furthermore, we have:
\begin{displaymath}
[\sigma_t^2|z_t,z_{t-1},\theta_t] \sim \mathcal{IG}(\alpha_t,\frac{Q_t}{2})
\end{displaymath}
where 
\begin{displaymath}
\alpha_t = \frac{n_t-p_t-q_{t-1}}{2}
\end{displaymath}
\begin{displaymath} Q_t = (z_t - H_t \hat{\lambda}_t)^T R_t^{-1}(D_t)(z_t - H_t \hat{\lambda}_t) \end{displaymath}

\begin{displaymath}
  \hat{\lambda}_t = \mathbb{E}[(\beta_{\rho_{t-1}},\beta_t)|z_t,z_{t-1},\theta_t,\sigma_t^2] 
\end{displaymath}
The REML estimator of $\sigma_t^2$ is $ \hat{\sigma}_t^2 = \frac{Q_t}{2\alpha_t }$ and we can estimate $\theta_t$ by minimizing the expression:
\begin{displaymath}
\mathrm{log}(|\mathrm{det}(R_t(D_t))|) + (n_t-p_t-q_{t-1})\mathrm{log}(\hat{\sigma}_t^2 )
\end{displaymath}

\subsection{Some important results about the covariance matrix $V_s$}\label{SomeImpV}

By sorting the experimental design sets as in Subsection \ref{ImpRes}, it can be shown that $\forall t = 2,\dots,s$ the inverse of the matrix $V_s$ has the form:
\begin{displaymath}\label{big_inverseA}
V_s^{-1} = \left(
\begin{array}{cc}
V_{s-1}^{-1} + \left(
\begin{array}{cc}
0 & 0 \\
0 &(\rho_{s-1}(D_s)\rho_{s-1}^T(D_s) ) \odot \frac{R_{s}^{-1}}{\sigma_{s}^2} \\
\end{array}
\right) &
- \left(
\begin{array}{c}
0 \\
(\rho_{s-1}(D_s)\mathbf{1}_{n_{s}}^T) \odot \frac{R_{s}^{-1}}{\sigma_{s}^2} \\
\end{array}
\right) 
\\
- \left(
\begin{array}{cc}
0 &
(\mathbf{1}_{n_{s}}\rho_{s-1}^T(D_s) )\odot \frac{R_{s}^{-1}}{\sigma_{s}^2} \\
\end{array}
\right) & 
\frac{R_{s}^{-1}}{\sigma_{s}^2} \\
\end{array}
\right)
\qquad
\end{displaymath}
with $V_1^{-1} = \frac{R_{1}^{-1}}{\sigma_{1}^2}$,  $V_{s-1}^{-1}$ an $(\sum_{i=1}^{s-1}{n_i} \times \sum_{i=1}^{s-1}{n_i})$ matrix and $R_{s}^{-1}$ an $(n_s \times n_s)$ matrix.
It can also be shown that:
\begin{displaymath}
t_s(x)^T V_s^{-1} = \left(
\begin{array}{c@{,}c}
\rho_{s-1}(x)t_{s-1}^T(x) V_{s-1}^{-1} -[0, \rho_{s-1}^T(D_s)\odot R_s(\{x\},D_s)R_s^{-1} ] &
R_s(\{x\},D_s)R_s^{-1}
\end{array}
\right)
\end{displaymath}

\subsection{Bayesian prediction for a code with 2 levels}

The equations for the Bayesian prediction when $\rho$ depends on $x$ can be directly derived from the Section \ref{Baypred2levels} by replacing $\rho$ with $\beta_\rho$ and noting that the design matrix $F$ is such that:

\begin{displaymath}
F=[F_\rho(D_2)\odot (z_1(D_2)\mathbf{1}_{p_\rho}^T) \quad F_2] 
\end{displaymath}

Finally, for the Bayesian prediction, we just have to adapt the integral (\ref{eq11}) :

\begin{displaymath}\label{Aeq11}
p(z_2(x)| z_1, z_2,\sigma_1^2,\sigma_2^2) = \int_{\mathbb{R}^{p_{\rho}+p_2}}{ p(z_2(x)|z_1, z_2,\beta_2,\beta_{\rho},\sigma_1^2,\sigma_2^2)p(\beta_{\rho},\beta_2|z_1, z_2,\sigma_2^2) \, d\beta_{\rho}d\beta_2 }
\end{displaymath}

\label{fin}

\end{document}